\documentclass[11pt,a4paper]{amsart}
\usepackage{amscd,amssymb,graphicx}
\author{Alice Fialowski}
\address{E\"otv\"os Lor\'and University\\
Budapest, Hungary} \email{fialowsk@cs.elte.hu}
\author{Michael Penkava}
\address{University of Wisconsin\\
Eau Claire, WI 54702-4004} \email{penkavmr@uwec.edu}
\thanks{The research of the authors was partially supported by grants
from the Mathematisches Forschungsinstitut Oberwolfach, OTKA T043641,
T043034 and the University of Wisconsin-Eau Claire.}

\newtheorem{thm}{Theorem}[section]

\theoremstyle{definition}
\newtheorem{dfn}[thm]{Definition}

\def \ph{\varphi}

\def \diag{\operatorname {diag}}

\def \ra{\rightarrow}

\def \hom{\mbox{\rm Hom}}
\def \ie{\hbox{\it i.e.}}

\def \tns{\otimes}

\def \mcom{,\cdots,}
\def \k{\mbox{$\mathbb K$}}

\def \C{\mbox{$\mathbb C$}}
\def \R{\mbox{$\mathbb R$}}
\def \Z{\mbox{$\mathbb Z$}}

\def\zt{\mbox{$\Z_2$}}

\def\sh{\operatorname{Sh}}

\def\inv{^{-1}}
\def\d{d}

\def\im{\operatorname{Im}}
\def\A{\mbox{$\mathcal A$}}

\def\linf{\mbox{$L_\infty$}}
\def\and{\mbox{ \rm and }}

\def\s#1{(-1)^{#1}}

\def\psa#1#2{\psi^{#1}_{#2}}

\def\inv{^{-1}}

\def\aut{\operatorname{Aut}}
\def\dinfty{\mbox{$d^\infty$}}
\def\dinf{\mbox{$d^\text{inf}$}}
\def\P{\mathbb P}

\begin{document}
\setlength{\multlinegap}{0pt}
\title[Deformations and Moduli Spaces of Lie Algebras]
{Formal Deformations, Contractions and Moduli Spaces of Lie algebras}%

\address{}%
\email{}%

\thanks{}%

\date{\today}
\begin{abstract}
Jump deformations and contractions of Lie algebras are inverse concepts, but the
approaches to their computations are quite different. In this paper, we contrast
the two approaches, showing how to compute the jump deformations from the miniversal
deformation of a Lie algebra, and thus arrive at the contractions. We also compute
contractions directly. We use the moduli spaces of real 3-dimensional and complex
3 and 4-dimensional Lie algebras as models for explaining a deformation theory
approach to computation of contractions.
\end{abstract}
\maketitle

\section{Introduction}
Deformations of analytic and algebraic objects is an old problem in
both mathematics and physics. In this paper we restrict ourselves to
the case of Lie algebras - one of the most important categories in
physics. On deformations theory of Lie algebras we refer to
\cite{fi1,fi3}.
The set of equivalence classes of Lie algebras over a fixed vector space is called
the moduli space of Lie algebras on that vector space.  In \cite{fp3,Ott-Pen}, the moduli
space of Lie algebras of dimension 3 was carefully analyzed, and in \cite{fp8},
a construction of the moduli space of 4-dimensional Lie algebras was given. The main
idea which we used in our analysis was the computation of the miniversal deformation,
which allows one to determine all possible deformations of the Lie
algebra (see \cite{fi2,ff2}). 

From the miniversal deformation, one can determine all jump deformations of a Lie
algebra. A jump deformation is precisely the inverse of a contraction of a Lie
algebra, so one can say that the miniversal deformations contain all the information
about contractions as well as other interesting information about the moduli space.

The point of view of deformation theory is a bit different from the point of view
of contractions.  When computing a contraction, one has a particular Lie algebra
in mind, and wants to know all Lie algebras which can jump to the one you have in
mind. This is quite different from the perspective of deformation theory, where one
is interested in seeing what the object of question deforms to.  Both perspectives
give valuable insights. We should refer to the recent work \cite{fm1,
fm2} which compare these two concepts and also give some examples. At the Workshop "Deformations and Contractions
in Mathematics and Physics" in Oberwolfach in January
2006, organized by Alice Fialowski, Marc de Montigny, Sergei Novikov
and Martin Schlichenmaier, researchers
from both mathematics and physics were brought together to share the ideas coming from
these two approaches. It was a valuable experience to the authors.

In this paper, we look at some of the examples we have previously studied from the
deformation point of view, and consider both the deformation and contraction points
of view.  We point out some of the advantages and disadvantages of these two approaches,
illustrating them with examples from moduli spaces of 3 and 4-dimensional Lie algebras.
We also give a description of the moduli space of real 3-dimensional Lie algebras. In
\cite{fp3, Ott-Pen, fp8}, only complex Lie algebras were studied.  We also
use a miniversal deformation approach to give a complete
description of the contractions of all 4-dimensional complex Lie algebras.
For analogous reasearch we refer to \cite{fm1} and \cite{fm2}.

Some contractions can be computed by use of diagonal matrices, and in fact, if one
is clever about a choice of a basis, one can always compute contractions in this manner
\cite{Weim}.  However, it is not true that one can compute all contractions given a
fixed basis by using diagonal matrices, and we give some counterexamples.
\section{Preliminaries}
Let $V$ be a \zt-graded vector space defined over $\C$, and denote
the even and odd parts of $V$ by $V_e$ and $V_o$, respectively.
The parity reversion $W=\Pi V$ is given by $W_e=V_o$ and
$W_o=V_e$; in other words, we reverse the parity of homogeneous
elements of $V$ to obtain $W$. Let $\pi:V\ra W$ be the identity;
note that it is an odd map.  Lie algebras are defined as
antisymmetric maps $V\tns V\ra V$; that is, as elements of
$\hom(\bigwedge^2(V),\ra V)$. These antisymmetric maps induce symmetric
maps $W\tns W\ra W$, so a Lie algebra determines an element in
$\hom(S^2(W),W)$, where $S^k(W)$ is the $k$-th symmetric power of $W$.

If we let
$E(V)=\bigoplus_{k=0}^\infty \bigwedge^k(V)$ be the \zt-graded
exterior algebra of $V$, then $E(V)$ can be identified
as a vector space, with
$S(W)$, where $S(W)=\bigoplus_{k=0}^\infty S^k(W)$ is the
\zt-graded symmetric algebra of $W$.  The algebra structures of these
spaces do not coincide, unless $V$ is a totally even space. For ordinary
Lie algebras, since $E(V)$ and $S(W)$ are isomorphic as algebras, it
is perfectly reasonable to work in the $E(V)$ picture, which is the point
of view in the classical literature, but for Lie superalgebras, there is
a big advantage in working with the $S(W)$ picture, so we will adopt this
point of view in this paper.

Any element of $\hom(E(V),V)$ can be identified with an element of $C(W)=\hom(S(W),W)$,
and therefore, a Lie algebra structure on $V$ determines an element $d$
in $\hom(S^2(W),W)$.  The structure $d$ satisfies the following relation, corresponding
to the Jacobi identity.
\begin{equation*}
d(d(a,b),c)+\s{bc}d(d(a,c),b)+\s{a(b+c)}d(d(b,c),a)=0,
\end{equation*}
where $\s{ab}$ is minus one to the power of the product of the parities of $a$ and $b$.
Let $C^k(W)=\hom(S^k(W),W)$. Then $C(W)=\prod_{k=0}^\infty C^k(W)$.
We define a product $\circ$ on $C(W)$ as follows. If $\ph\in C^k(W)$ and
$\psi\in\C^l(W)$, then $\ph\circ\psi\in C^{k+l-1}(W)$ is given by
\begin{equation*}
\ph\circ\psi(w_1\cdots w_n)=\sum_{\sigma\in\sh(l,k-1)}
\epsilon(\sigma)\ph(\psi(w_{\sigma(1)}\cdots w_{\sigma(l)})w_{\sigma(l+1)}\cdots w_{\sigma(n)}),
\end{equation*}
where $n=k+l-1$ and $\epsilon(\sigma)$ is a sign determined by the rule
\begin{equation*}
w_{\sigma(1)}\cdots w_{\sigma(n)}=\epsilon(\sigma)w_1\cdots w_n.
\end{equation*}
The bracket $[\ph,\psi]$ is defined by
\begin{equation*}
[\ph,\psi]=\ph\circ\psi-\s{\ph\psi}\psi\circ\ph.
\end{equation*}
This bracket equips the space of cochains
$C(W)$ with the structure of a \zt-graded Lie algebra. In fact, it is
well known that $C(W)$ is naturally isomorphic to the space of coderivations of
the symmetric coalgebra $S(W)$ of $W$, which is a \zt-graded Lie algebra, and the bracket
introduced above is just the bracket of coderivations \cite{sta4}.
In fact, any element $\ph\in C^k(W)$ extends to a codervation $\ph:S(W)\ra S(W)$, which
by the formula
\begin{equation*}
\ph(w_1\cdots w_n)=\sum_{\sigma\in\sh(k,n-k)}
\epsilon(\sigma)\ph(w_{\sigma(1)}\cdots w_{\sigma(k)})w_{\sigma(k+1)}\cdots w_n.
\end{equation*}
Accordingly, the cochains in $C(W)$ will sometimes be referred to as coderivations.

The Jacobi identity for $d$ is precisely the condition that $[d,d]=0$.
Moreover, $d$ is odd.  We call such an odd element of $C(W)$ a codifferential. In
fact, the map $D:C(W)\ra C(W)$, given by $D(\ph)=[d,\ph]$ is a differential on $C(W)$,
and is a derivation with respect to the bracket on $S(W)$. That is,
$D([\ph,\psi])=[D(\ph),\psi]+\s{\ph}[\ph,D(\psi)]$.  Note that $D$ is odd. We call
$D$ the coboundary operator induced by $d$, and the homology
$H(D)$ of this differential, defined by
\begin{equation*}
H(d)=\ker(d)/\im(d)
\end{equation*}
is called the cohomology of $d$. Since $D:C^k(W)\ra C^{k+1}(W)$, we can also define
the $n$-th cohomology group
$H^n(d)$ by $$H^n(d)=\ker(D:C^n(W)\ra C^{n+1}(W))/
\im(D:C^{n-1}(W)\ra C^n(W)).$$
Only the odd part of
$H^2(d)$ and the even part of $H^3(d)$ play a role in the theory of deformations of $d$.
For ordinary Lie algebras, $H^2(d)$ is a completely odd space and $H^3(d)$ is completely
even, because the parity of a cochain in $C^k(W)$ in this case depends only on $k$.

An \linf\ algebra is defined as a codifferential $d$ in $C(W)$, the only difference being
that we do not restrict $d$ to lie in $C^2(W)$. In fact, we can express
$d=d_1+\cdots$, where $d_i\in C^i(W)$. (Here, we do not allow a $d_0$ term.)  However,
in the case of \linf\ algebras, the cohomology $H(d)$ cannot be decomposed into
subgroups $H^n(d)$, and all of the cohomology plays a role in the deformation theory.
For this paper, we will restrict to Lie algebras, for simplicity, but the main
constructions extend to Lie superalgebras and \linf\ algebras as well.

The importance of cohomology to deformations is illustrated by the notion of an
infinitesimal deformation $d_t$ of $d$, which is given by
\begin{equation*}
d_t=d+t\delta,
\end{equation*}
where $t^2=0$.  The Jacobi identity $[d_t,d_t]=0$ reduces to the cocycle condition
$D(\delta)=0$. To understand why cohomology arises in the classification, we need to
introduce the notion of infinitesimal equivalence.

If $g$ is an automorphism of $W$; \ie, an invetertible linear map,
then $g$ extends uniquely to an automorphism of $S(W)$,
that is, an invertible map, compatible with the coalgebra
structure of $S(W)$. Moreover, $g$ acts on $C(W)$ by the rule
$g^*(\ph)=g\inv\ph g$.  We define two codifferentials $d'$ and $d$ to be equivalent,
and write $d\sim d'$, if there is an automorphism $g$ such that $d'=g^*(d)$.
An infinitesimal automorphism of
$W$ is a map $g_t=1+t\lambda$ where $\lambda:W\ra W$ is linear. If we
extend  $\lambda$ as a coderivation of $S(W)$, then we have $g_t=\exp(t\lambda)$.
We can thus identify $\lambda$ with
a cochain in $C(W)$. Evidently, $g_t\inv=1-t\lambda$.  We have
$g_t^*(\ph)=\ph+t[\ph,\lambda]$ for any $\ph\in C(W)$.

Now suppose that $d_t=d+t\delta$ and $d_t'=d+t\delta'$.  Then $d_t\sim d_t'$ precisely
when there is some cochain $\lambda$ such that $\delta'=\delta+D(\lambda)$, in other
words, when $\delta'$ and $\delta$ belong to the same cohomology class.  This is why
we say that the infinitesimal deformations are classified by the cohomology.

A formal deformation $d_t$ is given by
a formal power series of the form $$d_t=d+t\delta_1+t^2\delta_2+\cdots.$$
The Jacobi identity for $d_t$ is equivalent to the relations
\begin{equation*}
D(\delta_n)=-\tfrac12\sum_{k+l=n}[\delta_k,\delta_l]=0.
\end{equation*}
for $n=1,\dots$.  If the relations above hold for all $n<m$, then the right hand side
of the equation is a cocycle, but the fact that it is a coboundary is nontrivial. One
says that $d_t$ is an $m$-th order deformation if the relations hold for $n\le m$,
and that the $m$-th order deformation extends to an $(m+1)$-th order deformation
if there is some $\delta_{m+1}$ satisfying the relation above for $n=m+1$. A formal
deformation is an $m$-th order deformation for all $m$.  From these remarks, we
see that there is a relationship between cohomology and formal deformations, but
it is less straightforward than in the infinitesimal case.

One can also consider deformations, both infinitesimal and formal, in which more than
one parameter $t$ appears. Suppose that $H^2=\langle \delta^1\mcom \delta^n\rangle$.
Here, we identify $H^2$ with a subspace of $C^2(W)$, and the cochains $\delta^i$ are called
representative cocycles for a basis of the cohomology.
Consider the infinitesimal deformation $\dinf$, given by
\begin{equation*}
\dinf=d+t_i\delta^i.
\end{equation*}
This particular infinitesimal deformation is called the universal infinitesimal deformation.
It has the nice property that it uniquely
generates all infinitesimal deformations in the following
sense. If $d_t$ is any infinitesimal deformation given by the parameter $t$, then $d_t$ is
equivalent to a deformation of the form $d+t\delta$ where $\delta=a_i\delta^i$. Then
there is an obvious map from $k[t_1\mcom t_n]\ra k[t]$, such that $t_i\mapsto a_it$. With
this assignment, the universal infinitesimal deformation determines $d_t$ uniquely.

One can apply a similar construction for formal deformations to arrive at what is called
a miniversal deformation of $d$. Let $\{\delta^i\}$ be a (pre)basis of $H^2$ as in the
infinitesimal case, and $\{\gamma^i\}$ be a prebasis of the 3-coboundaries. Then there
is a deformation of the form
\begin{align*}
\dinfty=d+t_i\delta^i +s_i\gamma^i,
\end{align*}
where $s_i$ is a formal power series in the variables $t_i$, with all terms in $s_i$ of order
at least 2, such that
\begin{equation*}
[\dinfty,\dinfty]=r_i\alpha^i,
\end{equation*}
where $\{\alpha^i\}$ is a basis of $H^3$ and $r_i$ are formal
power series in the $t_i$ starting with
terms of degree at least 2. The series $r_i$ are called the relations on the base of the
miniversal deformation.  If
these series converge in some neighborhood
of $0\in C^n$ and $(t_1\mcom t_n)\in C^n$ satisfies all the relations $r_i$, then
the deformation $d_t$ is a well defined codifferential.

In the examples we present,
it turns out that the relations on the base
are actually given by rational functions of the parameters,
which are defined at $t_i=0$, so that there is such a radius of convergence.
Thus, we are able to obtain a notion of neighborhood
of a codifferential $d$, consisting of those deformations arising from substituting
small values of the parameters into the miniversal deformation.

This notion of neighborhood does not give rise to a Hausdorff topology on the space of
equivalence classes of codifferentials on $W$, which is called the moduli space of
Lie algebras on $W$. It is possible for there
to be a formal deformation $d_t$, which is well defined as a codifferential for small
values of $t$, for which $d_0=d$, but for which $d_t\sim d'$ for all values of $t$ except
$t=0$, where $d'$ is not equivalent to $d$. This kind of a deformation family is called
a jump deformation. If a codifferential has a jump deformation to another codifferential,
then it is not a closed point in the moduli space. The existence of such points is why
the moduli space is not Hausdorff. This suggests that if we could somehow exclude the
jump deformations, we could introduce a more reasonable topological decomposition of
the moduli space of Lie algebras.

A deformation family $d_t$, where $d_t$ runs along a family
of nonequivalent codifferentials, is called a smooth deformation family. In this
article, we shall describe a stratification of the orbifold by smooth orbifolds,
where the smooth neighborhoods of a codifferential are given by smooth deformations.
The jump deformations provide a type of gluing operation between the strata. All
of the non-Hausdorff behavior of the moduli space is thus represented by the jump
deformations.

If $d_t$ is a jump deformation from $d$ to $d'$, so that $d_0=d$ and $d_t\sim d'$ for
$t\ne 0$, and $d'_t$ is a jump deformation from $d'$ to $d''$, there is no way to
compose these jump deformations directly to obtain a jump deformation $d''_t$ from
$d$ to $d''$.  Nevertheless, there is always such a jump deformation.  Thus jump
deformations are transitive.  Similarly, if $d_t$ is a jump deformation from
$d$ to $d'$, and $d_t'$ is a smooth family of deformations of $d'$, then there is
a smooth family $d_t''$ of deformations of $d$, and a smooth function $f(t)$ with
$f(0)=0$, such that $d_{f(t)}'=d_t''$. In this case, we say that the deformation
$d_t''$ factors through the jump deformation from $d$ to $d'$.

The stratification of the moduli space is obtained by considering as neighborhoods of
a codifferential $d$ those codifferentials which can be obtained as smooth deformations
of $d$ which do not factor through jump deformations.
The smooth deformations which factor
through a jump deformation always belong to a different stratum.  Jump deformations
are always a one way phenomenon. If $d$ jumps to $d'$, then it never happens that
$d'$ will jump to $d$. In fact, if $\dinfty(t)$ is
 the miniversal deformation of $d$, parameterized by $t=(t_1\mcom t_n)$
 then for a small enough neighborhood of
$0$ in the $t$ space, $\dinfty(t)$ is not equivalent to $d$ unless $t=0$. Thus a smooth
deformation $d_t$ of $d$, given by a curve arising from the miniversal deformation,
cannot have $d_t\sim d$, for small nonzero values of $t$. This explains why jump deformations
are one way.

There is an obvious way in which projective geometry enters the description of the moduli
space of Lie algebras on $W$, because $ad\sim d$, for any nonzero $a$.  However, in our
study of 3 and 4 dimensional complex Lie algebras
\cite{fp3,fp8}, we discovered that the orbifolds which form
the strata of the moduli space have a natural structure of a projective orbifold. The
projective structure which arises is not simply the consequence of the identification
of codifferentials with their multiples; it is a more subtle relationship.

The goal of this paper is to illustrate the phenomena we have described with some
examples, and to describe the moduli space of 3 and 4 dimensional Lie algebras. These
moduli spaces were studied in detail in \cite{fp3,fp8,Ott-Pen}.  When constructing
moduli spaces of Lie and \linf\ algebras, we had noticed that jump deformations were
always transitive. In January 2006, at a workshop in Oberwolfach on Deformations and
Contractions in Mathematics and Physics, we learned that the physics notion of contraction
is equivalent to the mathematical notion of jump deformation. That is to say, a contraction
and a jump deformation are inverse notions.  If a codifferential $d$ has a jump deformation
to $d'$, then $d'$ contracts to the codifferential $d$.  Using this equivalence, and a
result of E. Weimar-Woods \cite{Weim}, we will give a simple proof of the
transitivity of jump deformations. However, the transitivity of jump deformations is
equivalent to the transitivity of contractions, which was pointed out in \cite{Weim}.

A common approach to constructing formal deformations is to use Massey products, and one
can solve the problem of whether a particular infinitesimal deformation extends to a formal
deformation using this approach.  The advantage of studying the miniversal deformation
is that it gives all of the formal deformations at once, and makes it possible to
analyze the structure of the moduli space locally. The construction of the miniversal
deformation which we use here first appeared in \cite{ff2},
and was extended to infinity
algebras in \cite{fp1}.
\section{Moduli spaces of Complex Lie algebras}
Let us consider the moduli space of complex
three dimensional Lie algebras. This space is quite
simple in structure, consisting of a one-parameter family of solvable Lie algebras,
and three special Lie algebras.

To describe this moduli space, let us introduce some notation.
If $W=\langle w_1\mcom w_n\rangle$ is an $n$-dimensional odd vector space,
$1\le i\le n$ and
$I=(i_1\mcom i_k)$, where $1\le i_1<\cdots < i_k\le n$, then $\ph^I_i$ denotes
the element of $C^k(W)$ given by $$\ph^I_i(w_J)=\delta^I_J w_i,$$ where
$w_J=w_{j_1}\cdots w_{j_k}$. When $k$ is odd, $\ph^I_i$ is even. Similarly, when
$k$ is even $\ph^I_i$ is odd, and to emphasize the difference between even and odd
elements, we shall denote it by $\psi^I_i$ instead.  Elements in $C^2(W)$ are all
odd, because we assume that $W$ is odd. This picture corresponds to ordinary Lie
algebras, because then the vector space $V$ is even.

If $W=\langle w_1,w_2,w_3\rangle$ is a 3-dimensional completely
odd vector space, then $S^2(W)=\langle
w_1w_2,w_1w_3, w_2w_3\rangle$, and $C^2(W)$ is thus 9-dimensional. In terms
of this basis, an element $\ph\in C^2(W)$, given by
\begin{align*}
\ph&=a_{11}\psa{12}1+a_{12}\psa{13}1+a_{13}\psa{23}1\\
&+a_{21}\psa{12}2+a_{22}\psa{13}2+a_{23}\psa{23}2\\
&+a_{31}\psa{12}3+a_{32}\psa{13}3+a_{33}\psa{23}3
\end{align*}
is given by the $3\times3$ matrix $A=(a_{ij})$.

The classification of three-dimensional algebras is classical,
for example, it appears in \cite{jac,ov}. In order to give the correct
stratification of the moduli space, it is necessary to realign the classical
decomposition slightly.  In \cite{fp8}, we gave the special points the
names $d_1$, $d_2$ and $d_3$, and the elements of the family were denoted by
$d_2(\lambda:\mu)$, where $(\lambda:\mu)$ are projective coordinates.  Let us
identify these elements with the classical notation.  There is only one nontrivial
nilpotent Lie algebra (up to isomorphism), which is $\mathfrak n_3(\C)$.
This algebra coincides with our $d_1$. There
is a family $\mathfrak r_{3,\alpha}(\C)$ of solvable Lie algebras, which coincides
with $d_2(\lambda:\mu)$ for $\alpha=\mu/\lambda$, except when $\alpha=1$.
The Lie algebra $\mathfrak r_{3,1}(\C)$ coincides with our special point $d_2$, while
the solvable Lie algebra $\mathfrak r_3(\C)$ corresponds to our point $d_2(1:1)$.
The Lie algebra $\mathfrak r_{3,0}(\C)$ is also denoted as $\mathfrak r_2(\C)\oplus\C$,
where $\mathfrak r_2(\C)$ is the nontrivial 2-dimensional Lie algebra.  Finally,
the simple Lie algebra $\mathfrak{sl}_2(\C)$ coincides with our $d_3$.

The real difference
between our classification and the usual one is that we interchange the elements
$\mathfrak r_3$ and the elements $\mathfrak r_{3,1}$. This interchange arises from
the necessity of aligning the elements in the moduli space into strata that are
distinguished by jump deformations. The element $d_2(1:1)$ belongs in the strata
with the family $d_2(\lambda:\mu)$, rather than the element $d_2$, because there
is a jump deformation from $d_2$ to $d(1:1)$, instead of the other way around. Both
of them have smooth deformations along the family $d_2(\lambda:\mu)$, but the smooth
deformations of $d_2$ along this family factor through the jump deformation to
$d_2(1:1)$.

Actually, the first hint that the family might be misaligned can be seen
in the behaviour of the cohomology of the Lie algebras.  The dimension of
$H^2(\mathfrak r_{3,1})$ is 3, while the dimension of $H^2(\mathfrak r_{3,\alpha})$ is
generically equal to 1. The dimension of $H^2(\mathfrak r_3)$ is 1, which is
appropriate for an element in the family.

Another anomaly in the cohomology occurs for
$\mathfrak r_{3,-1}$, when the dimension of $H^2$ jumps to 2. In the classical picture,
one considers only the elements $\mathfrak r_{3,\alpha}$ where $|\alpha|\le 1$. More
precisely, we find that $\mathfrak r_{3\alpha}\sim\mathfrak r_{3,1/\alpha}$. In our
notation, $d_2(\lambda:\mu)\sim d_2(\mu:\lambda)$.
Since $(\lambda:\mu)$ parameterizes the Riemann sphere $\P^1(\C)$, we obtain an
action of the symmetric group $\Sigma_2$ on this Riemann sphere, with the equivalence
classes of the codifferentials $d_2(\lambda:\mu)$ parameterized by the orbifold
$\P^1(\C)/\Sigma_2$. Note that we have obtained the stratum as a projective orbifold.

There are precisely two orbifold points in $\P^1(\C)/\Sigma_2$,
the points $(1:1)$ and $(1:-1)$.
(Recall that the orbifold points are the points whose stabilizer is nontrivial.)
It is thus not surprising that something special should occur at the orbifold points.
The codifferential $d_2(1:1)$ is special, because there is a jump deformation from the
codifferential $d_2$ to it, while the codifferential $d_2(1:-1)$ is special for the
reason that it has a jump deformation to the codifferential $d_3$.  In fact, this
jump deformation is well known to physicists, because it corresponds to a contraction
of the simple Lie algebra $\mathfrak{sl}_2(\C)$.

The nilpotent Lie algebra $d_1$ has jump deformations to every codifferential in the
family $d_2(\lambda:\mu)$, as well as to the simple Lie algebra $d_3$. In fact, the
jump deformation from $d_1$ to $d_3$ is an example of a transitive jump, because
it factors through the jump deformation from $d_1$ to $d_2(1:-1)$. It is not surprising
that the nilpotent Lie algebra should deform to the solvable Lie algebras. Nilpotent
Lie algebras are the least rigid in terms of their deformations, while simple Lie
algebras are completely rigid.

Consider the solvable Lie algebra $d_2$, which is represented by the codifferential
$d_2=\psa{13}1 +\psa{23}2$. Explicitly, in terms of the basis $\{w_1,w_2,w_3\}$, the
Lie algebra structure, in standard bracket notation, is given by
\begin{equation*}
[w_1,w_3]=w_1,\quad [w_2,w_3]=w_2,
\end{equation*}
with all other brackets vanishing.
In \cite{fp3} the miniversal deformation of
$d_2$ was calculated. Its formula, as given in \cite{Ott-Pen}, is
\begin{equation*}
\d_2^\infty=\psa{13}1(1+t_1)+\psa{23}2+\psa{13}2t_2+\psa{23}1t_3.
\end{equation*}
In terms of standard bracket notation, the deformed algebra is given by the bracket rules
\begin{equation*}
[w_1,w_3]=(1+t_1)w_1+t_2w_2,\quad [w_2,w_3]=t_3w_1+w_2.
\end{equation*}
It is convenient to express the miniversal deformation of $d_2$ by the matrix
$A=\left[\begin{smallmatrix}0&1+t_1&t_3\\0&t_2&1\\0&0&0\end{smallmatrix}\right]$.

The 3 parameters in the versal deformation arise from the fact that
$H^2(d_2)$ is 3-dimensional, which
means that the tangent space of the versal deformation is 3-dimensional. This is an
interesting situation, because the moduli space of 3-dimensional Lie algebras consists
of a 1-dimensional piece, and three 0-dimensional pieces, so it is a bit disconcerting
to find that the dimension of the tangent space is larger than the dimension of the
moduli space.

A partial explanation for this phenomenon is as follows. When classifying infinitesimal
deformations, one considers the action of the group of infinitesimal automorphisms of
$W$, that is, the maps $g=1_W+t\lambda$, where $\lambda\in\hom(W,W)$. Under the
action of this group, the infinitesimal deformations of the form $d_t=d+t\delta$, are
classified by the cohomology $H^2(d)$. However, the automorphism group $\aut(d)$,
consisting of the automorphisms $g$ of $W$ such that $g^*(d)=d$ acts on $H^2(d)$,
and thus on the set of infinitesimal deformations,
and it is isomorphism classes under this action
which really classify the nonequivalent deformations of $d$. Thus the tangent space
should really be considered as an orbifold, in terms of the action of this group.

If one considers which of the codifferentials $d_2^\infty$ is equivalent to, one
finds that
\begin{equation*}
d_2^\infty\sim d_2(\alpha:\beta),
\end{equation*}
where
\begin{equation*}
\alpha,\beta=1+1/2\,t_{{1}}\pm1/2\,\sqrt {{t_{{1}}}^{2}+4\,t_{{3}}t_{{2}}}
\end{equation*}
Thus, the deformation is equivalent to $d_2(1:1)$ precisely when $t_1^2+4t_2t_3=0$,
and not all the parameters vanish.
For example, the 1-parameter family of deformations
$$d_t=d_2^\infty(0,t,0)=\psa{13}1+\psa{23}2+\psa{13}2t$$ is equivalent to
$d_2(1:1)$ whenever $t\ne0$. This is an example of a jump deformation. If $g_t$ is the
automorphism of $W$ given by the matrix
\begin{equation*}
\left[ \begin {array}{ccc}
0&1&1\\
\noalign{\medskip}t&0&1\\
\noalign{\medskip}0&0&1
\end {array} \right],
\end{equation*}
then $g_t^*(d_t)=d_2(1:1)$. To compute this, we first compute the matrix $Q_t$ representing
$g_t:S^2(W)\ra S^2(W)$.
We have
\begin{align*}
g_t(w_1w_2)&=tw_2w_1=-tw_1w_2\\
g_t(w_1w_3)&=tw_2(w_1+w_2+w_3)=-tw_1w_2+tw_2w_3\\
g_t(w_2w_3)&=w_1(w_1+w_2+w_3)=w_1w_2+w_1w_3,
\end{align*}
so that
\begin{equation*}
Q_t=\left[\begin{matrix}
-t&-t&1\\0&0&1\\0&t&0
\end{matrix}\right].
\end{equation*}
If $A_t$ represents the matrix for $d_t$ and $A'$ the matrix for $d_2(1:1)$, then
\begin{equation*}
A_t=\left[\begin{matrix}
0&1&0\\0&t&1\\0&0&0
\end{matrix}\right]
\qquad
A'=\left[\begin{matrix}
0&1&1\\0&0&1\\0&0&0
\end{matrix}\right].
\end{equation*}
In matrix form the condition $g_t^*(d_t)=d_2(1:1)$ is simply $G_t\inv A_t Q_t=A'$,
which is easily verified.

Turning this process around, we obtain $(g\inv)^*(d(1:1))=d_t$. This formula represents
the fact that $d(1:1)$ contracts to the Lie algebra $d_2$. Let us recall the definition
of a contraction.
\begin{dfn}
Let $g_t$ be a family of automorphisms of $V$, defined in a punctured neighborhood of
zero.  If $\lim_{t\ra0}(g_t^*(d'))\sim d$ exists and
$d$ is not equivalent to $d'$, then the Lie algebra given by $d$ is said to be a contraction of
the Lie algebra given by $d'$.
\end{dfn}

Since $\lim_{t\ra0}(d_t)=d_2$, we see that $d_2$ is a contraction of the Lie algebra
$d(1:1)$. Note that the jump deformation $d_t$ from $d_2$ to $d(1:1)$ is by no means
unique. In fact, if one considers any curve $\gamma(t)$ on the surface $t_1^2+4t_2t_3=0$
satisfying $\gamma(0)=0$, then the deformation $d_t=d_2^\infty(\gamma(t))$ represents
a jump deformation from $d_2$ to $d(1:1)$, so there is a corresponding contraction
of $d(1:1)$ to $d_2$.

To determine all possible contractions of a codifferential $d'$ by finding all
automorphisms $g_t$ such that $\lim_{t\ra0}(g_t^*(d'))=d$ exists can be a
daunting task.  However, what is interesting is not the method of obtaining
the contraction, but simply the contracted object. The following theorem \cite{Weim},
due to E. Weimar-Woods, makes the task of computing the contractions much simpler in
practice.
\begin{thm}[Weimar-Woods, 2000]\label{weimar}
If there is a contraction from $d'$ to $d$, where $d$ and $'d$ are Lie algebra
structures on a finite dimensional space $W$, then there is a basis of $W$ and
an automorphism $g_t$ of $W$, which has matrix $\diag(t^{\lambda_1}\mcom t^{\lambda_n})$,
where the $\lambda_i$ are integers, such that $g^*(d')\sim d$.
\end{thm}
This theorem makes it possible to determine all Lie algebras $d$ which arise as contractions
of $d'$, even if the classification of Lie algebra structures on $W$ is not known. Let us
illustrate this idea by considering the contractions of $d'=d(1:1)$. The matrix of
$d'$ is $\left[\begin{smallmatrix}0&1&1\\0&0&1\\0&0&0\end{smallmatrix}\right]$. Suppose
that the matrix $G$ of $g_t$ is given by $G=\diag(t^a,t^b,t^c)$. Then the matrix
$Q$ of $g:S^2(W)\ra S^2(W)$ is $Q=\diag(t^{a+b},t^{a+c},t^{b+c})$. Thus the matrix of
$A'=G\inv A Q$ is
$\left[\begin{smallmatrix}0&t^c&t^{b+c-a}\\0&0&t^c\\0&0&0\end{smallmatrix}\right]$.
in order that $\lim_{t\ra0}g^*(d')$ exists, we must have $c\ge0$ and $b+c-a\ge0$.
Clearly, if both of the inequalities are strict, this describes the uninteresting
contraction to the zero codifferential. Also, if both inequalities are equalities,
this does not describe a contraction, since the original codifferential is not changed.
Thus we have two nontrivial contractions, to the codifferential $d_1$,
given by the matrix $\left[\begin{smallmatrix}0&0&1\\0&0&0\\0&0&0\end{smallmatrix}\right]$
 and
to the codifferential $d_2$, given by the matrix $\left[\begin{smallmatrix}
0&1&0\\0&0&1\\0&0&0\end{smallmatrix}\right]$.

As is illustrated by the example, given a finite dimensional Lie algebra, there are only
a finite number of Lie algebras which can arise as contractions of the Lie algebra. The
converse statement about jump deformations is not true.  In fact, the Lie algebra $d_1$
has jump deformations to every 3-dimensional Lie algebra except for $d_2$.  It is also
said that the multiplication in a contracted Lie algebra is ``more abelian''. More
precisely, one can say that the cohomology of a contracted Lie algebra is higher
dimensional, and that the contracted Lie algebra has ``more deformations'' than the
original Lie algebra. To show this, we will analyze the miniversal deformation of
a Lie algebra more carefully.

Suppose that $\{\delta^1\mcom \delta^m\}$ is a pre-basis of $H^2(d)$, in other words,
we assume that $\delta^i$ are 2-cocycles whose images in $H^2(d)$ are a basis, and that
$\{\gamma^1\mcom \gamma^n\}$ is a pre-basis of the 3-coboundaries, so that the
$D(\gamma^i)$ give a basis for $B^3(W)=D(C^2(W))$. Then the miniversal deformation
can be given in the form
\begin{equation*}
d^\infty=d+t_i\delta^i+x_j\gamma^j,
\end{equation*}
where $x_j$ are formal power series in the parameters $t_i$, whose lowest order terms
are of degree 2; \ie,
\begin{equation*}
x_i=a^{jk}_i t_jt_k+ a^{jkl}_it_jt_kt_l+\cdots.
\end{equation*}
If $\{\alpha^i\}$ is a pre-basis of $H^3(W)$ and $\{\tau^i\}$ is a prebasis of the
$B^4(W)$, then
\begin{equation*}
[\dinfty,\dinfty]=r_i\alpha^i+u_i\tau^i,
\end{equation*}
where the relations $r_i$ are formal power series in the parameters $t_i$, whose lowest
order terms are of degree 2, and the $u_i$ are formal power series in the $t_i$ which
are contained in the ideal in $\k[[t_1\mcom t_m]]$ generated by the relations. The
 ring $\A=\k[[t_1\mcom t_m]]/(r_i)$ is called the base of the miniversal deformation.
The existence of a miniversal deformation of this form for Lie algebras
was proved in \cite{ff2}, and for infinity algebras in \cite{fp1}.

A formal 1-parameter deformation $d_t$ of $d$ is given by any homomorphism
$f:\A\ra k[[t]]$, such that $f(t_i)$ is a formal power series in $t$ with
nonzero constant term, and is defined by
\begin{equation*}
d_t=d+a_i(t)\delta^i+b_i(t)\gamma^i,
\end{equation*}
where $a_i(t)=f(t_i)$ and $b_i(t)=f(x_i)$ are formal power series in $t$.
If the power series $f(t_i)$ and $f(x_i)$ converge in a neighborhood of zero,
then the deformation is said to be analytic.
In principle, the computation of the miniversal deformation might be
expected to be a fairly intractable problem.

One can proceed to compute the miniversal deformation order by order, and hope
that the process terminates after a finite number of steps. In many of the examples
which the authors have studied, this procedure does work.  However, another idea
is to write the formula for the versal deformation $\dinfty$ as above, with unknown
coefficients $x_i$. Then one computes that
\begin{equation*}
[\dinfty,\dinfty]=r_i\alpha^i+s_i\beta^i+u_i\tau^i,
\end{equation*}
where the $\beta^i$ are a basis of the 3-coboundaries.  The deformation will be miniversal
if $s_i=0$ for all $i$. If we have chosen $\beta^i=D(\gamma_i)$, then we see that
\begin{equation*}
s_k=2x_k+a^{ij}_kt_it_j +b^{ij}_kt_ix_j+c^{ij}_kx_ix_j.
\end{equation*}
The number of equations is exactly equal to the number of variables $x_i$. From the form
of the equations, there is a solution near $x=t=0$.  However, these equations are quadratic
in the $x_i$, so it is not clear that it is possible to solve them in any systematic manner.
Surprisingly, for any three or four dimensional example, it turns out that there is not
only a solution, but the solution expresses the $x_i$ as rational functions of the $t_i$.
It is not known to us whether this property is true in general, but it is true for every
example which we have constructed.

After solving for the $x_i$, one substitutes their values into the expressions above for
the $r_i$, to obtain the relations as functions of the parameters $t_i$.
The fact that the $u_i$ are equal to zero mod the $r_i$ follows from the
construction in \cite{ff2}.

Now let us suppose that $d_t'$ is an analytic deformation of $d'$. then $d_t'=d'+t\ph(t)$.
Suppose that $d=\lim_{t\ra0}g^*(d')$ is a contraction of $d'$. Then $d_t=g_t^*(d')$
is a jump deformation of $d$. We do not know that $\lim_{t\ra0}g_t^*(d_t')$ exists.
However, let us suppose that $g$ is expressed as a diagonal matrix in integer powers
of $t$.  Then, if $k$ is a large enough odd positive integer exponent, $g^*(t^k\ph(t^k))$ will
be given by positive powers of $t$ only, and therefore, its limit as $t\ra0$ will be zero.
Thus, the deformation $\tilde d_t=g^*_t(d_{t^k}')$ is a well defined deformation of $d$.

Note that if $d'_t$ is a jump deformation from $d'$ to $d''$, then $\tilde d_t$ is a
jump deformation of $d$ to $d''$, which shows the transitivity of jump deformations.
Actually, the transitivity of contractions is known as well \cite{Weim}, and is even
more obvious.  When $d_t'$ is a smooth deformation, that is, when the codifferentials
$d_t$ are not isomorphic as $t$ varies, then $\tilde d_t$ is also a smooth deformation,
and $d_{t^k}'\sim \tilde d_t$. We say that $\tilde d_t$ factors through the jump deformation
$d'_t$. We summarize this analysis in the theorem below.
\begin{thm}
Suppose that $d'_t$ is a deformation of $d'$ and that there is a jump deformation of
$d$ to $d'$.  Then, for a sufficiently large positive integer $k$, there is a deformation $d_t$
of $d$ such that $d_t\sim d'_{t^k}$.
\end{thm}

As an example, let us consider the case $d'=d(1:1)$ with deformation
$d_t'$ represented by the matrix
$A=\left[\begin{smallmatrix}0&1&1\\0&t&1\\0&0&0\end{smallmatrix}\right]$.
Let $G=\diag(t^a,t^b,1)$ represent $g_t$, where $b-a>0$. Since $g_t^*(d')$
is given by the matrix
$\left[\begin{smallmatrix}0&1&t^{b-a}\\0&0&1\\0&0&0\end{smallmatrix}\right]$,
$g_t$ determines a contraction of $d'$ to $d=d_2$. The matrix of
$g_t^*(d'_{t^k})$ is
$\left[\begin{smallmatrix}0&1&t^{b-a}\\0&t^{k-(b-a)}&1\\0&0&0\end{smallmatrix}\right]$.
It follows that for $k\ge b-a$, the deformation $d_t=g_t^*(d'_{t^k})$ is a well defined
deformation of $d_2$.  Note that in this case, if we choose $b-a=1$, then we can set
$k=1$. However, it is not clear from our argument whether one can always find an appropriate
$g_t$ so that we can set $k=1$. The key issue is that if $d''$ is a deformation of $d'$,
then it is also a deformation of $d$. This motivates the following definition.
\begin{dfn}
Suppose that $d_t\sim d'_t$ in some punctured neighborhood of $t=0$, $d_0=d$ $d'_0=d'$ and there
is a jump deformation from $d$ to $d'$.  Then we say that the deformation $d_t$ factors
through the jump deformation.
\end{dfn}

For the sake of completeness, we will now give a complete description of the miniversal
deformations of the three dimensional Lie algebras, and compute their contractions.
If
$A=(a_{ij})$ is an arbitrary matrix of a codifferential, and $G=\diag(t^a,t^b,t^c)$
is a diagonal matrix of an automorphism $g_t$ of $\C^3$, then $A'=G\inv A Q$, where $Q$ is the
matrix of $g_t:S^2(W)\ra S^2(W)$, has the form
\begin{equation*}\label{contractform}
A'=\left[
\begin{array}{rrr}
t^ba_{11}&t^ca_{12}&t^{b+c-a}a_{13}\\
t^aa_{21}&t^{a+c-b}a_{22}&t^ca_{23}\\
t^{a+b-c}a_{31}&t^aa_{32}&t^ba_{33}
\end{array}
\right].
\end{equation*}
In order to obtain a nontrivial contraction, all of the powers of $t$
corresponding to nonzero entries in $A$
must be nonnegative, at least one, but not all, of the powers must be zero. To obtain the
matrix of the contracted codifferential, one simply lets $t=0$ in $A'$.

Of course, even
though two resulting matrices may be different, the codifferentials may still be equivalent,
so it is necessary to check this.  However, note that it is not necessary to know the complete
classification of the Lie algebras in order to check whether different matrices give rise
to equivalent contractions.  This is one of the strengths of the contraction method, because
it can be applied to determine all nonequivalent contractions of a codifferential, without
a knowledge of the complete classification of the moduli space.

\subsection{The simple Lie algebra $\mathfrak{sl}_2(\C)$}
The simple Lie algebra $\mathfrak{sl}_2(\C)$ is represented
by the codifferential $d_3=\psa{12}3+\psa{13}2+\psa{23}1$, with
matrix $\left[\begin{smallmatrix}0&0&1\\0&1&0\\1&0&0\end{smallmatrix}\right]$.
 As the cohomology of this
codifferential vanishes completely, the
versal deformation is simply $d_3^\infty=d_3$ which is not interesting. On the other hand,
one computes immediately that the matrices
$\left[\begin{smallmatrix}0&0&0\\0&1&0\\1&0&0\end{smallmatrix}\right]$,
$\left[\begin{smallmatrix}0&0&1\\0&0&0\\1&0&0\end{smallmatrix}\right]$,
and
$\left[\begin{smallmatrix}0&0&1\\0&1&0\\0&0&0\end{smallmatrix}\right]$, all arise from
contractions.  These matrices arise from codifferentials which are equivalent to
$d_2(1:-1)$.
\subsection{The solvable Lie algebra $\mathfrak{r}_{3,1}(\C)$}
This is given by the codifferential $d_2=\psa{13}1+\psa{23}2$. We have already discussed
the versal deformation of this codifferential. Note that from the form \eqref{contractform}
of a diagonal contraction, it follows that there are no nontrivial contractions
of $d_2$.
\subsection{The solvable Lie algebra $\mathfrak{r}_{3,-1}(\C)$}
This is given by the codifferential $d(1:-1)=\psa{13}1+\psa{23}1-\psa{23}2$. This
codifferential is unique in the family $d(\lambda:\mu)$ in that its cohomology $H^2$ is
two dimensional, which means that its versal deformation is given by a two parameter
family $$d^\infty(1:-1)=\psa{13}1(1+t_1)+\psa{23}1-\psa{23}2+\psa{12}3t_2,$$
whose matrix is
$\left[\begin{smallmatrix}0&1+t_1&1\\0&0&-1\\t_2&0&0\end{smallmatrix}\right]$.
There is one relation on the base: $t_1t_2=0$. This means that either $t_1=0$ or
$t_2=0$. In the former case, the versal deformation is equivalent to the simple
Lie algebra $d_3$, whenever $t_2\ne0$, so this gives a jump deformation. When
$t_2=0$, the versal deformation is equivalent to the codifferential $d_2(1+t_1:-1)$, which
means that as $t$ changes, the deformation moves along the family $d_2(\lambda:\mu)$.

For contractions, we note that one can obtain the matrices
$\left[\begin{smallmatrix}0&0&1\\0&0&0\\0&0&0\end{smallmatrix}\right]$, corresponding
to the codifferential $d_1$, and
$\left[\begin{smallmatrix}0&1&0\\0&0&-1\\0&0&0\end{smallmatrix}\right]$, which is equivalent
to the codifferential $d_2(1:-1)$ again. Thus, only the first matrix yields a
contraction.
\subsection{The solvable Lie algebras $\mathfrak{r}_{3,r}(\C)$, $\mathfrak{r}_3(\C)$, and
$\mathfrak{r}_2(\C)\oplus\C$}
Recall that the codifferential $d_2(\lambda:\mu)$ represents the Lie algebra
$\mathfrak{r}_{3,\mu/\lambda}$ when $\lambda\ne0$ and $\lambda\ne\mu$. When $\lambda=0$,
this represents the Lie algebra $\mathfrak{r}_2(\C)\oplus\C$, and when $\lambda=\mu$, this
represents the codifferential $\mathfrak{r}_3(\C)$. Except for the case
when $\lambda=\mu$, the versal deformation is given by
$d^\infty(\lambda:\mu)=\psa{13}1(\lambda+t_1)+\psa{23}1+\psa{23}2\mu$, with matrix
$\left[\begin{smallmatrix}0&\lambda+t_1&1\\0&0&\mu\\0&0&0\end{smallmatrix}\right]$.
It is very clear from the form of the matrix that the versal deformation is isomorphic
to $d_2(\lambda+t_1:\mu)$, so the deformations move along the same family. This is what
determines the neighborhood structure of elements in the family. Note that since there
are no jump deformations, these smooth deformations do not factor through a jump deformation,
so it is natural to identify the neighborhoods of $d_2(\lambda:\mu)$ as being given
by the family.

When $\lambda=\mu$, it turns out that the cocycle
$\psa{13}1$ is a coboundary, so cannot be used as in the generic case to
parameterize the versal deformation. It is strange that the deformation
$d_t=d(1:1)+\psa{13}1t$ varies smoothly along the family, although its leading term is
a coboundary.

It is easy to see that the only nontrivial contractions of $d_2(\lambda:\mu)$, when
$\lambda\ne\mu$ are to $d_1$, representing the Lie algebra $\mathfrak{n}_3$. However,
when $\lambda=\mu$, the matrix
$\left[\begin{smallmatrix}0&1&0\\0&0&1\\0&0&0\end{smallmatrix}\right]$ arises by the
contraction process, and this matrix corresponds to the codifferential $d_2$, which therefore
is a contraction of $d_2(1:1)$.
\subsection{The nilpotent Lie algebra $\mathfrak{n}_3(\C)$}
The nilpotent Lie algebra $\mathfrak{n}_3(\C)$ is
represented by the codifferential $d_1=\psa{23}1$. As the dimensional
of $H^2(d_1)$ is 5, it is not surprising that $d_1$ has a lot of deformations;
in fact, it deforms to
every 3-dimensional Lie algebra except $d_2$. The versal deformation is given by the matrix
$\left[\begin{smallmatrix}0&0&1\\-t_2&t_5&t_3\\t_4&t_2&t_1\end{smallmatrix}\right]$,
and there are two relations on the base:
$t_1t_5-t_2t_3=0$ and $ t_1t_2+t_3t_4=0$.
We will not give explicitly formulas for all the jump deformations. To determine them,
one first solves the relations explicitly. Then, for a solution of the relations, one
determines what codifferential is represented by the corresponding matrix. All of this
is easy to do using a computer algebra system. Note that some of the deformations are
not jump deformations, but run along the family $d_2(\lambda:\mu)$. These are examples
of smooth deformations which factor through a jump deformation.
There are no nontrivial contractions of $d_1$.
\section{Real 3-dimensional Lie algebras}

In the classification of  complex 3-dimensional Lie algebras, one can proceed as follows. Either
the algebra is simple, in which case it is isomorphic to $\mathfrak{sl}_2(\C)$, which
is represented by the codifferential $d_3$, or it is solvable, so it is an extension
of the 1-dimensional Lie algebra by a 2-dimensional one.

It turns out that one only needs
to consider the case of an extension of a 1-dimensional Lie
algebra  by the abelian 2-dimensional Lie algebra, because the extensions by the nontrivial
2-dimensional Lie algebra do not give any additional nonequivalent codifferentials.
The matrix of the extension can be given in the form
$A=\left[\begin{smallmatrix}0&A'\\0&0\end{smallmatrix}\right]$, where $A'$ is an arbitrary
$2\times2$ matrix. Two such extensions are equivalent precisely when the matrices are similar,
up to multiplication by a constant.  As a consequence, we can reduce everything to the
Jordan decomposition of the matrix.

In fact, the codifferential
$d_2(\lambda:\mu)$ is given by the matrix
$A'=\left[\begin{smallmatrix}\lambda&1\\0&\mu\end{smallmatrix}\right]$,
$d_2$ is given by the identity matrix, and $d_1$ is given by the matrix
$A'=\left[\begin{smallmatrix}0&1\\0&0\end{smallmatrix}\right]$.  The interpretation of
$d_2(\lambda:\mu)$ as a $\P^1(\C)/\Sigma_2$ is a consequence of fact that equivalence
is given by similarity.

The same pattern can be observed in higher dimensions.  For an $n+2$-dimensional Lie
algebra, there is a stratum that is given as the orbifold $\P^n/\Sigma_{n+1}$, where
the action of $\Sigma_{n+1}$ is given by permuting the projective coordinates in
$\P^n$.

For real Lie algebras, one has to make the following modifications of the theory. First,
there are two nonisomorphic simple Lie algebras
$\mathfrak{sl}_2(\R)$, represented by the
codifferential $d_3=\psa{12}3+\psa{13}2+\psa{23}1$ and $\mathfrak{su}_2$,
represented by the codifferential $d'_3=\psa{12}3-\psa{13}2+\psa{23}1$.

Secondly, while it is true that any extension of $\R$ by a 2-dimensional real
Lie algebra is equivalent to one given by an extension by an abelian Lie algebra, so
that it is determined by the similarity class of a $2\times 2$ matrix $A'$,
the rational canonical form determines the similarity classes of real matrices.
Codifferentials of the form $d_2(\lambda:\mu)$, determined by
the matrices $\left[\begin{smallmatrix}\lambda&1\\0&\mu\end{smallmatrix}\right]$
are still nonequivalent over $\C$. These codifferentials are parameterized
by $\P^1(\R)/\Sigma_2$. However, there are problems which arise in using this family
as part of the description of the moduli space.

It turns out to be more effective to work with the matrices from the
rational canonical form, which can be expressed in the form
$A_{x,y}=\left[\begin{smallmatrix}0&1\\x&y\end{smallmatrix}\right]$.
It is easy to check that
the codifferentials corresponding to $A_{x,y}$ and $A_{t^2x,ty}$ are equivalent, for any
$t$. To get a single family of codifferentials, parameterized by the action of $\Sigma_2$
on $\P^1$, it is convenient to give the family $d(\lambda:\mu)$ as follows:
\begin{equation*}
d(\lambda:\mu)=
\begin{cases}
\psa{13}2\lambda+\psa{23}1\lambda+2\psa{23}2\mu&\text{if $\lambda\ge 0$}\\
-\psa{13}2\lambda+\psa{23}1\lambda+2\psa{23}2\mu&\text{if $\lambda<0$}.
\end{cases}
\end{equation*}
It is easy to check that that $d(\lambda:\mu)\sim d(t\lambda:t\mu)$ for $t>0$.
Note that our space is not really projective, corresponding to the quotient of
$\R^2-\{0\}$ by $\R_{+}$, rather than $\R_*$.
Furthermore $d(\lambda:\mu)\sim d(\lambda:-\mu)$, so we obtain
an action of $\Sigma_2$ on the parameter space, determining the codifferentials
up to equivalence.  There are two orbifold points in this action, $(1:0)$ and $(-1:0)$.

A justification for the seemingly artificial gluing together of two types of codifferential
at the point $d(0:1)$ is given by studying the versal deformation of the Lie algebra, which
is
\begin{equation*}
d^\infty(0:1)=\psa{13}1t+\psa{23}2.
\end{equation*}
The versal deformation is equivalent to the codifferential $d(\lambda:\mu)$
where $\mu=1+t$, and $\lambda=\sqrt{-t}$ if $t<0$, and $\lambda=-\sqrt{t}$ if $t>0$.
Thus the two pieces of $d(\lambda:\mu)$ are glued together at $d(0:1)$ by means
of the versal deformation.

The points $d(1:0)$ and $d(-1:0)$ both correspond to the same point $d_2(1:-1)$
in the complex case, so it not surprising that each of them has a jump deformation
to a simple Lie algebra. In fact, the contractions of the 3-dimensional
simple Lie algebras were computed in \cite{conatser} and a complete list of the contractions
is given in \cite{ov}.

Let us consider the real algebra $\mathfrak{sl}_2(\R)$, which is given by the
same codifferential $d_3=\psa{12}3+\psa{13}2+\psa{23}1$ as the complex
algebra $\mathfrak{sl}_2(\C)$, with matrix
$\left[\begin{smallmatrix}0&0&1\\0&1&0\\1&0&0\end{smallmatrix}\right]$.
We have already computed the contractions of $d_3$, previously,
but now we need to identify the real algebras associated to them.
The complex algebras are all isomorphic to $d_2(1:-1)$, but there are two real
versions of this algebra, $d(1:0)$ and $d(0:1)$.
The matrix $A$ contracts to the matrices
$\left[\begin{smallmatrix}0&0&1\\0&1&0\\0&0&0\end{smallmatrix}\right]$,
$\left[\begin{smallmatrix}0&0&1\\0&0&0\\1&0&0\end{smallmatrix}\right]$, and
$\left[\begin{smallmatrix}0&0&0\\0&1&0\\1&0&0\end{smallmatrix}\right]$.
The first matrix is just the matrix of $d(1:0)$. The second matrix gives a
codifferential equivalent to $d(-1:0)$, while the third is equivalent to $d(1:0)$ again.
Thus, there are two distinct contractions of $\mathfrak{sl}_2(\R)$.

The real algebra $\mathfrak{su}_2$ is given by the codifferential $d_3'
=\psa{12}3-\psa{13}2+\psa{23}1$, with matrix
$\left[\begin{smallmatrix}0&0&1\\0&-1&0\\1&0&0\end{smallmatrix}\right]$. Its contractions
are given by the matrices
$\left[\begin{smallmatrix}0&0&1\\0&-1&0\\0&0&0\end{smallmatrix}\right]$,
$\left[\begin{smallmatrix}0&0&1\\0&0&0\\1&0&0\end{smallmatrix}\right]$, and
$\left[\begin{smallmatrix}0&0&0\\0&-1&0\\1&0&0\end{smallmatrix}\right]$. Clearly,
the first one is just the matrix of $d(-1:0)$. The other two are also equivalent
to $d(-1:0)$. As a consequence, there is just one contraction of $\mathfrak{su}_2$.

Note that $d(-1:1)\sim d_2(1:1)$. The matrix of $d(-1:1)$ is
$\left[\begin{smallmatrix}0&0&1\\0&-1&2\\0&0&0\end{smallmatrix}\right]$. Since a contraction
by a diagonal matrix has the effect of setting some of the coefficients in the matrix to zero,
it is impossible to obtain a matrix which is equivalent to $d_2$ by using a diagonal matrix
to perform the contraction, with respect to this basis. In fact, if we consider the
automorphism $g_t$, given by the matrix $G_t=
\left[\begin{smallmatrix}1/t&0&0\\1/t&1&0\\0&0&1\end{smallmatrix}\right]$, then
$g_t^*(d(-1:1)$ is given by the matrix
$\left[\begin{smallmatrix}0&1&t\\0&0&1\\0&0&0\end{smallmatrix}\right]$. Therefore,
$\lim_{t\ra0}g_t^*(d(1:1))\sim d_2$. This example illustrates an important limitation
of the contraction method in \cite{Weim}, that one has to be clever about finding a
basis in which the matrix of the automorphism producing the contraction is diagonal.
As far as we can see, this limitation is a serious one, because the advantage of
computation of contractions by use of diagonal matrices is in the ease of computation,
but if the procedure misses some of the contractions, it is inadequate to solving the
problem posed in \cite{Weim5}, that of finding a simple class of contractions which produce
all possible contractions.
\section{Complex 4-dimensional algebras}
In \cite{fp8}, the moduli space of 4-dimensional complex Lie algebras was studied in
detail, and a decomposition into strata consisting of orbifolds, connected by jump
deformations was given. Miniversal deformations for the Lie algebras were computed, so
all contractions of these Lie algebras can be read off from the jump deformations.
We will use the basis $\{w_1w_2,w_1w_3,w_2w_3,w_1w_4,w_2w_4,w_3w_4\}$ for $S^2(W)$,
where $W=\langle w_1,w_2,w_3,w_4\rangle$ is a completely odd space of dimension 4.
The moduli space of Lie algebras can be decomposed into
one 2-dimensional orbifold, two 1-dimensional orbifolds, and 6 special points.
The decomposition is as follows.
\begin{enumerate}
\item $\mathbf{d_3(\lambda:\mu:\nu)}$: with matrix $\left[\begin{smallmatrix}
0&0&0&{\it \lambda}&1&0\\\noalign{\medskip}0&0&0&0&{\it \mu}&1\\\noalign{\medskip}0&0&0&0&0&{
\it \nu}\\\noalign{\medskip}0&0&0&0&0&0
\end{smallmatrix}\right]$ is a 2-dimensional family of codifferentials, where
$(\lambda:\mu:\nu)$ are projective coordinates, and the action of the group $\Sigma_3$,
by permuting the coordinates gives equivalent codifferentials. Thus this family is
parameterized by the orbifold $\P^1(\C)/\Sigma_3$.\vspace{.1in}
\item $\mathbf{d_1(\lambda:\mu)}$: given by
$\left[ \begin {smallmatrix}
0&0&1&{\it \mu}+{\it \lambda}&0&0\\\noalign{\medskip}0&0&0&0&{\it
\lambda}&1\\\noalign{\medskip}0&0&0&0&0
&{\it \mu}\\\noalign{\medskip}0&0&0&0&0&0
\end {smallmatrix} \right]
$ is a 1-dimensional family of codifferentials, given by projective coordinates
$(\lambda:\mu)$, with an action of $\Sigma/2$, given by permutation of coordinates.
Thus this family is parameterized by $\P^1(\C)/\Sigma_2$.\vspace{.1in}
\item $\mathbf{d_3(\lambda:\mu)}$: given by the matrix
$\left[\begin{smallmatrix}
0&0&0&{\it \lambda}&0&0\\\noalign{\medskip}0&0&0&0&{\it \lambda}&1\\\noalign{\medskip}0&0&0&0&0
&{\it \mu}\\\noalign{\medskip}0&0&0&0&0&0
\end{smallmatrix}\right]$ is the other 1-dimensional family, given by projective
coordinates $(\lambda:\mu)$. This family does not have an action of $\Sigma_2$, so
it is parameterized simply by $\P^1(\C)$.\vspace{.1in}
\item $\mathbf{d_1}$: given by
$\left[\begin{smallmatrix}
 0&0&0&0&1&0\\\noalign{\medskip}0&0&0&0&0&0
 \\\noalign{\medskip}0&0&0&0&0&0\\\noalign{\medskip}0&0&0&0&0&0
\end{smallmatrix}\right]$ is the nilpotent Lie algebra $\mathfrak{n}_3(\C)\oplus\C$.
\vspace{.1in}
\item $\mathbf{d_1^\sharp}$: given by
$\left[\begin{smallmatrix}
0&0&1&2&0&0\\\noalign{\medskip}0&0&0&0&1&0
\\\noalign{\medskip}0&0&0&0&0&1\\\noalign{\medskip}0&0&0&0&0&0
\end{smallmatrix}\right]
$ is a solvable algebra.\vspace{.1in}
\item $\mathbf{d_2^\star}$: given by
$\left[\begin{smallmatrix}
 0&0&0&0&1&0\\\noalign{\medskip}0&0&0&0&0&1
 \\\noalign{\medskip}0&0&0&0&0&0\\\noalign{\medskip}0&0&0&0&0&0
 \end{smallmatrix}\right]
$ is the nilpotent Lie algebra $\mathfrak{n}_4(\C)$.\vspace{.1in}
\item $\mathbf{d_2^\sharp}$: given by
$\left[\begin{smallmatrix}
 1&0&0&0&0&0\\\noalign{\medskip}0&0&0&0&0&0
 \\\noalign{\medskip}0&0&0&0&0&1\\\noalign{\medskip}0&0&0&0&0&0
\end{smallmatrix}\right]$ is the Lie algebra $\mathfrak r_2(\C)\oplus\mathfrak r_2(\C)$.
\vspace{.1in}
\item
$\mathbf{d_3}$: given by
$\left[\begin{smallmatrix}
0&0&1&0&0&0\\\noalign{\medskip}0&1&0&0&0&0
\\\noalign{\medskip}1&0&0&0&0&0\\\noalign{\medskip}0&0&0&0&0&0
\end{smallmatrix}\right]$ is the Lie algebra $\mathfrak{sl}_2(\C)\oplus\C$.\vspace{.1in}
\item
$\mathbf{d_3^*}$: given by
$\left[\begin{smallmatrix}
0&0&0&1&0&0\\\noalign{\medskip}0&0&0&0&1&0
\\\noalign{\medskip}0&0&0&0&0&1\\\noalign{\medskip}0&0&0&0&0&0
\end{smallmatrix}\right]$ is another solvable Lie algebra.
\end{enumerate}
Complete details about the miniversal deformations of these algebras, and an explanation
for the decomposition is given in \cite{fp8}, and we do
not reproduce this information here. Our goal here is to study the contractions of
the codifferentials.

\section{Contractions of $\mathfrak{sl}_2(\C)\oplus\C$.}
The codifferential $d_3$ is the same codifferential as in the 3-dimensional case, but in
four dimensions, it picks up additional contractions.
Note that $d_3$ is still rigid as a 4-dimensional Lie algebra, although
its cohomology does not vanish completely. The following series of jump
deformations give the complete contraction picture for $d_3$.
\begin{equation*}
d_1\rightsquigarrow d_2^*\rightsquigarrow d_3(1:-1:0)\rightsquigarrow d_1(1:-1)\rightsquigarrow
d_3.
\end{equation*}
Every 3-dimensional Lie algebra determines a 4-dimensional Lie algebra in a trivial
way, so we may consider the 3-dimensional Lie algebras as part of the 4-dimensional
moduli space.

Note that the $d_1$ in the 4-dimensional list is equivalent to the $d_1$ in the
3-dimensional case, and $d_3(\lambda:\mu:0)$ is equivalent to the Lie algebra given
by $d_2(\lambda:\mu)$, so two of the jumps on this list are already known from the
3-dimensional picture, and they account for all of the diagonal contractions of
$d_3$, in terms of the usual basis.
To understand the other contractions a bit better,
let us analyze the contraction to $d_1(1:-1)$. The the matrix of the miniversal deformation
of $d_1(1:-1)$ is
$\left[\begin{smallmatrix}
 0&0&1&t_{{2}}&0&0\\\noalign{\medskip}t_{{1}}&t_{{1}}-1/2\,t_{{1}}t_{{2}}&0&0&1+t_{{2}}&1\\\noalign{\medskip}0
&-t_{{1}}&0&0&0&-1\\\noalign{\medskip}0&0&0&t_{{1}}t_{{2}}&0&0
\end{smallmatrix}\right]$.
The relation on the base is $t_1t_2=0$. If we set $t_2=0$ and $t_1=t$ we obtain a
deformation $d_t$ given by
$\left[\begin{smallmatrix}
 0&0&1&0&0&0\\\noalign{\medskip}t&t&0&0&1&1\\\noalign{\medskip}0&-t&0&0&0&-1
 \\\noalign{\medskip}0&0&0&0&0&0
\end{smallmatrix}\right]
$.
If we let $g_t$ be given by  $G_t=
\left[\begin{smallmatrix}
 t&0&0&-1\\\noalign{\medskip}0&\sqrt {t}&0&0
 \\\noalign{\medskip}0&\sqrt {t}&\sqrt {t}&0\\\noalign{\medskip}0&0&0
&1
\end{smallmatrix}\right]$,
then $d_t=g_t^*(d_3)$, so $d_3$ contracts to $d_1(1:-1)$ using this transformation.
The square roots in the formula for $G_t$ are not important, because we can replace
$\sqrt t$ with $t$ in this matrix and obtain a contraction using $d_{t^2}$, which
clearly has the same limit.

Now, $d_3$ is equivalent to the codifferential obtained
by substituting $t=1$ in $d_t$. It is clear that we obtain $d_1(1:-1)$ as a contraction
using the basis in which $d_3$ has this matrix, but it is not obvious from the original
expression of $d_3$ that one should consider such a basis.  We determined this basis
from the jump deformations of $d_1(1:-1)$, but what we would have hoped to be
able to do is determine the contractions of $d_3$ without this knowledge.

\subsection{Contractions of $\mathfrak{r}_2(\C)\oplus\mathfrak{r}_2(\C)$}
The following series of jump deformations gives the complete picture of the contractions
of the codifferential $d_2^\sharp$, representing the Lie algebra
$\mathfrak{r}_2(\C)\oplus\mathfrak{r}_2(\C)$.
\begin{equation*}
d_1\rightsquigarrow d_2^\star\rightsquigarrow
d_3(1:0:1)\rightsquigarrow
d_1(1:0)\rightsquigarrow d_2^\sharp.
\end{equation*}
The algebra $d_2^\sharp$ is the only completely rigid 4-dimensional Lie algebra, which
means that its cohomology vanishes in all dimensions. Rigid Lie algebras have a complex
contraction picture.
\subsection{Contractions of $d_1(1:1)$}
Here the contractions do not come from a single line of jump deformations. We have
\begin{align*}
&d_1\rightsquigarrow d_2^\star\rightsquigarrow d_3(1:1:2) \rightsquigarrow d_1(1:1)\\
&d_1\rightsquigarrow d_1^\sharp\rightsquigarrow d_1(1:1).
\end{align*}
The first line is a special case of the generic contraction picture for codifferentials
of the form $d_1(\lambda:\mu)$, but the second line represents line of contractions which
only apply to this special point.

The codifferential $d_1(1:1)$ is analogous to the 3-dimensional
codifferential $d_2(1:1)$, in the way it behaves in its family. As a member of the
family $d_1(\lambda:\mu)$, it has no special properties with
respect to deformations, but there is a element
outside the family which has a jump deformation to it.  It is the element $d_1(1:-1)$
in the family $d_1(\lambda:\mu)$
which has extra deformations.  Note that this behaviour parallels the situation
with $d_2(\lambda:\mu)$, because the points $d_1(1:1)$ and $d_1(1:-1)$ are the
orbifold points in this family, and they are the ones where unusual behaviour occurs.

\subsection{Contractions of $d_1(\lambda:\mu)$}
Codifferentials of the form $d_1(\lambda:\mu)$ do not have many deformations, but they all
arise as jump deformations from the family $d_3(\lambda:\mu:\lambda+\mu)$. With the
exception of the codifferential $d_1(1:-1)$, which has
a jump deformation to $d_3$, all other deformations of $d_1(\lambda:\mu)$ simply move
along this family.
\begin{equation*}
d_1\rightsquigarrow d_2^\star \rightsquigarrow d_3(\lambda:\mu:\lambda+\mu)\rightsquigarrow
d_1(\lambda:\mu)
\end{equation*}
Note that in the action of $\Sigma_3$ on $\P^2$, the subgroup $\Sigma_2$ consisting of
the permutations of the first two coordinates preserves the $\P^1$ given by
$(\lambda:\mu:\lambda+\mu)$. As a consequence, $d_3(\lambda:\mu:\lambda+\mu)$ is parameterized
by $\P^1/\Sigma_2$, the same parameterization of $d_1(\lambda:\mu)$.

\subsection{Contractions of $d_3(\lambda:\lambda:\mu)$}
This special case of the codifferentials of type $d_3(\lambda:\mu:\nu)$ has a more
complex contraction picture than usual for this type.
deformations that run along the family $d_3(\lambda:\mu:\nu)$.
We have
\begin{equation*}
d_1\rightsquigarrow
d_2^\star\rightsquigarrow
d_3(\lambda:\mu)\rightsquigarrow
d_3(\lambda:\lambda:\mu)
\end{equation*}

These points $(\lambda:\lambda:\mu)$ are
orbifold points on $\P^2$, and the atypical contraction picture of the points
of the form
$d_3(\lambda:\lambda:\mu)$  resembles the behaviour
of the orbifold point $d_2(1:1)$ in the 3-dimensional case.
One does not observe
anything special about the deformations of codifferentials of the type
$d_3(\lambda:\lambda:\mu)$.

\subsection{Contractions of $d_3(\lambda:\mu:\nu)$}
Generically, there are not many contractions of $d_3(\lambda:\mu:\nu)$.
We have
\begin{equation*}
d_1\rightsquigarrow d_2^\star \rightsquigarrow d_3(\lambda:\mu:\nu)
\end{equation*}
Certain subfamilies of the family $d_3(\lambda:\mu:\nu)$ have extra jump deformations
or extra contractions, but generically, the deformations are only along the family.
\subsection{Contractions of $d_3(1:1)$}
In general, $d_3(\lambda:\mu)$ has no contractions. However, we obtain one special case.
\begin{equation*}
d_3^\star\rightsquigarrow d_3(1:1).
\end{equation*}
On the other hand $d_3(\lambda:\mu)$ has jump deformations to the family
$d_3(\lambda:\mu:\nu)$ along the subfamily $d_3(\lambda:\lambda:\mu)$. Note that
the action of $\Sigma_3$ on $\P^2$ identifies three copies
of $\P^1$, given by the action on $(\lambda:\lambda:\mu)$, which mutually intersect
only in the point $(1:1:1)$, so that  $d_3(\lambda:\lambda:\mu)$ has
no net nontrivial action of $\Sigma_3$. As a consequence, it is not surprising that
the parameter space of $d_3(\lambda:\mu)$
also has no group action, and is thus simply $\P^1$.
\section{A Complicated Versal Deformation}
To illustrate some of the difficulties with the approach to computing contractions
using the miniversal deformations of the objects in the moduli space, we give an
example of a 4-dimensional codifferential whose versal deformation is quite
complicated. For the codifferential $d_1=\psa{24}1$, the dimension of $H^2$ is 13,
while the dimension of $H^3$ is 10.  This means that there will be 10 relations
on the base of the versal deformation, which will be a 13-parameter algebra.

The matrix of the versal deformation is
$$\left[\begin{matrix}
0&t_{{6}}&0&t_{{4}}&1&0\\\noalign{\medskip}-t_{{11}}+t_{{4}}t_{{2}}&-t_{{9}}t_{{12}}+t_{{3}}t
_{{1}}+t_{{3}}t_{{2}}-t_{{9}}t_{{6}}&t_{{6}}+t_{{12}}&t_{{10}}&0&t_{{3
}}\\\noalign{\medskip}t_{{8}}&t_{{9}}t_{{2}}&t_{{1}}&t_{{13}}&0&t_{{9}
}\\\noalign{\medskip}t_{{7}}&-t_{{9}}t_{{5}}+t_{{12}}t_{{1}}+t_{{2}}t_
{{6}}+t_{{2}}t_{{12}}&t_{{5}}&t_{{11}}&t_{{2}}&t_{{12}}
\end{matrix}\right]$$
This complexity of this matrix is not the main source of the difficulty.  What makes it much
more difficult to work with is that there are relations on the base of the versal deformation.
The parameters need to satisfy 10 different relations, and finding the solutions to these
relations leads to a lot of complications.

When we solved these relations using Maple, we
came up with 48 distinct solutions. Each solution to the relations needs to be substituted
in the matrix, and then one has to consider which codifferential the solution is equivalent
to.  With a large number of parameters, determining the equivalence class of such a matrix
is not an easy task, even for the computer.  For the case of $d_1$, we were able to determine
all of the nonequivalent deformations, including the jump deformations, so we can determine
from this information which codifferentials contract to $d_1$. However, as the reader
can imagine, this example is about at the limit of complexity which can be handled by
current software.

The 10 relations on the base are
\begin{align*}
0=&
t_{{13}}t_{{1}}-t_{{13}}t_{{2}}-t_{{8
}}t_{{4}}+t_{{9}}t_{{8}}\\
0=&
-t_{{2}}t_{{11}}+{t_{{2}}}^{2}t_{{4}}-t_{{7}}
t_{{4}}+t_{{13}}t_{{5}}+t_{{12}}t_{{8}}\\0=&
-t_{{12}}t_{{4}}-2\,t_{{9}}t_{
{6}}+t_{{3}}t_{{1}}+t_{{3}}t_{{2}}-t_{{9}}t_{{12}}\\0=&
-t_{{1}}t_{{6}}-t_{{5}}t_{{4}}+t_
{{12}}t_{{1}}+t_{{2}}t_{{6}}+t_{{2}}t_{{12}}-t_{{9}}t_{{5}}\\0=&
-t_{{1}}t_{{11}}+
t_{{1}}t_{{4}}t_{{2}}-t_{{13}}t_{{5}}-2\,t_{{8}}t_{{6}}-t_{{12}}t_{{8}
}-t_{{9}}t_{{7}}\\
0=&-t_{{10}}t_{{2}}+t_{{8}}t_{{3}}+
t_{{4}}t_{{11}}-{t_{{4}}}^{2}t_{{2}}+t_{{13}}t_{{6}}+t_{{13}}t_{{12}}\\0=&
-t_{{
10}}t_{{1}}-t_{{8}}t_{{3}}-t_{{9}}t_{{2}}t_{{4}}+t_{{13}}t_{{6}}+t_{{
11}}t_{{9}}-t_{{13}}t_{{12}}\\0=&
t_{{1}}t_{{9}}t_{{5}}-t_{{12}}{t_{{1}}}^{2}-t_{{1}}t_
{{2}}t_{{6}}+{t_{{2}}}^{2}t_{{6}}+{t_{{2}}}^{2}t_{{12}}-2\,t_{{11}}t_{
{5}}\\&+t_{{5}}t_{{4}}t_{{2}}-2\,t_{{7}}t_{{6}}-2\,t_{{12}}t_{{7}}\\
0=&
2\,t_{{11}}t_{{3}}-
2\,t_{{3}}t_{{4}}t_{{2}}+t_{{4}}t_{{9}}t_{{12}}-t_{{4}}t_{{3}}t_{{1}}+
t_{{4}}t_{{9}}t_{{6}}-2\,t_{{12}}t_{{10}}\\&+{t_{{9}}}^{2}t_{{12}}-t_{{9}
}t_{{3}}t_{{1}}+{t_{{9}}}^{2}t_{{6}}\\0=&
-t_{{9}}t_
{{2}}t_{{12}}+t_{{2}}t_{{3}}t_{{1}}+t_{{3}}{t_{{2}}}^{2}-2\,t_{{2}}t_{
{9}}t_{{6}}-t_{{7}}t_{{3}}+t_{{4}}t_{{9}}t_{{5}}-t_{{4}}t_{{12}}t_{{1}
}\\&-t_{{4}}t_{{2}}t_{{6}}-t_{{4}}t_{{2}}t_{{12}}-t_{{10}}t_{{5}}+t_{{11}
}t_{{6}}+{t_{{9}}}^{2}t_{{5}}-t_{{9}}t_{{12}}t_{{1}}
\end{align*}
Four of the relations above contain only quadratic terms, while the rest
contain some cubic terms as well.  Because there are no higher order terms,
the versal deformation could have been computed by computing the deformation
order by order up to the third order. Many of our examples had relations which
were rational in the parameters, and in those cases, one could not calculate
the versal deformation order by order.
\section{Conclusions}
Although jump deformations and contractions are inverse concepts, the approaches
to their computation are quite different.  Each of these approaches has some
advantages and disadvantages.

The method of computation of jump deformations by computing miniversal
deformations is guaranteed to determine every object which contracts to the
object for which the miniversal deformation is calculated.  In this sense,
the computation of miniversal deformations contains all of the deformation
information, including all information about contractions.  However, it is not
so easy to compute the miniversal deformation except for algebras of low dimension.
To determine all the jump deformations which are contained in the
miniversal deformation is not easy either, mainly because it is not very simple
to determine when a family of deformations are equivalent to each other. The
difficulties that arise in using the miniversal deformation approach are mainly
due to computational complexity.

There is no satisfactory general method of computing all contractions
directly. Contractions
using diagonal matrices are simple to compute, and allows one to determine some of
the contractions easily.  However, as we have seen in this paper, there are cases
where a diagonal matrix is not sufficient to compute all of the contractions, given
a specific choice of basis of the underlying space, and it is not clear how to choose
a basis that will yield all of them by this method. On the other hand, while a Lie
algebra may have jump deformations to an infinite number of nonequivalent algebras,
a finite dimensional Lie algebra has only a finite number of contractions. Therefore,
by experimenting with different bases, one may reasonably expect to find
them all. Many different approaches, besides the use of diagonal matrices have
been tried, and they have led to very successful results.  The main difficulty with
the direct approach to computing contractions is that there is no general method
of determining them all.

The point of view is an important issue as well. For example, if one is interested in
computing all the contractions of the Poincar\'e algebra, then miniversal deformations
are not the right approach, because you don't know which algebras to compute the
miniversal deformations of. There is no classification of 10 dimensional Complex
Lie algebras, so it is impossible to proceed in this manner. On the other hand, if
one is interested in computing the deformations of the Galilean algebra, in particular,
to show that it deforms into the Poincar\'e algebra, as was done in \cite{fig},
a computation of the versal deformation is very useful. In fact, in the article
\cite{fig}, the versal deformation was not quite computed, but many deformations
of the Galilean algebra were computed using a partial computation of this
versal deformation.

Given the complexity of finding both deformations and contractions, we imagine that
both methods will continue to be valuable in computations. In higher dimensions, there
is no classification of the Lie algebras, so any special method of finding deformations
or contractions of particular Lie algebras may be illuminating.  It is possible that
in the future, a simple general method of computing all contractions of a Lie algebra
may be discovered.  Advances in computer hardware will make the computation of versal
deformations easier in the future as well.

\bibliographystyle{amsplain}
\providecommand{\bysame}{\leavevmode\hbox to3em{\hrulefill}\thinspace}
\providecommand{\MR}{\relax\ifhmode\unskip\space\fi MR }
\providecommand{\MRhref}[2]{%
  \href{http://www.ams.org/mathscinet-getitem?mr=#1}{#2}
}
\providecommand{\href}[2]{#2}

\end{document}